\documentclass[lettersize,journal]{IEEEtran}
\usepackage{amsmath,amsfonts}
\usepackage{algorithmic}
\usepackage{algorithm}
\usepackage{array}
\usepackage[caption=false,font=normalsize,labelfont=sf,textfont=sf]{subfig}
\usepackage{textcomp}
\usepackage{stfloats}
\usepackage{url}
\usepackage{verbatim}
\usepackage{graphicx}
\usepackage{subfig}
\usepackage{cite}
\usepackage{amsmath}

\DeclareMathOperator*{\argmin}{arg\,min}

\begin{document}

\title{Performance Metrics for Multi-Objective Optimisation Algorithms Under Noise}

\author{Juergen Branke
\thanks{University of Warwick, UK, juergen.branke@wbs.ac.uk}
}


\maketitle

\begin{abstract}
This paper discusses the challenge when evaluating multi-objective optimisation algorithms under noise, and argues that decision maker preferences need to be taken into account. It demonstrates that commonly used performance metrics are problematic when there is noise, and proposes two alternative performance metrics that capture selection error of the decision maker due to mis-estimating a solution's true fitness values.
\end{abstract}


\section{Introduction}
\IEEEPARstart{B}{eing} able to accurately measure the performance of an algorithm is fundamental to algorithm comparison and thus to scientific progress in algorithm design. For single-objective problems, this is relatively straightforward, one can simply compare the performances of the solutions returned by the algorithms. For multi-objective problems this is more challenging, as the algorithm returns a set of Pareto-optimal solutions, with different trade-offs regarding the underlying objectives. Researchers have developed a range of performance metrics for deterministic multi-objective problems, with the Hypervolume (HV), the Inverted Generational Distance (IGD) and the R2 metric being the most popular ones. These are all unary metrics which aggregate the performance of a set of solutions in a single value, capturing aspects such as closeness to the Pareto frontier, spread, and uniformity of the distribution of solutions along the Pareto frontier.

Measuring performance in a noisy single-objective problem is straightforward - one can simply use the true (noiseless) fitness of the solution returned by the algorithm. This is readily available for most artificial benchmark problems and can be estimated to arbitrary precision by averaging over a sufficiently large number of evaluations otherwise.

\begin{figure}
    \centering
    \subfloat[Error by exclusion]{
    \includegraphics[width=0.3\textwidth]{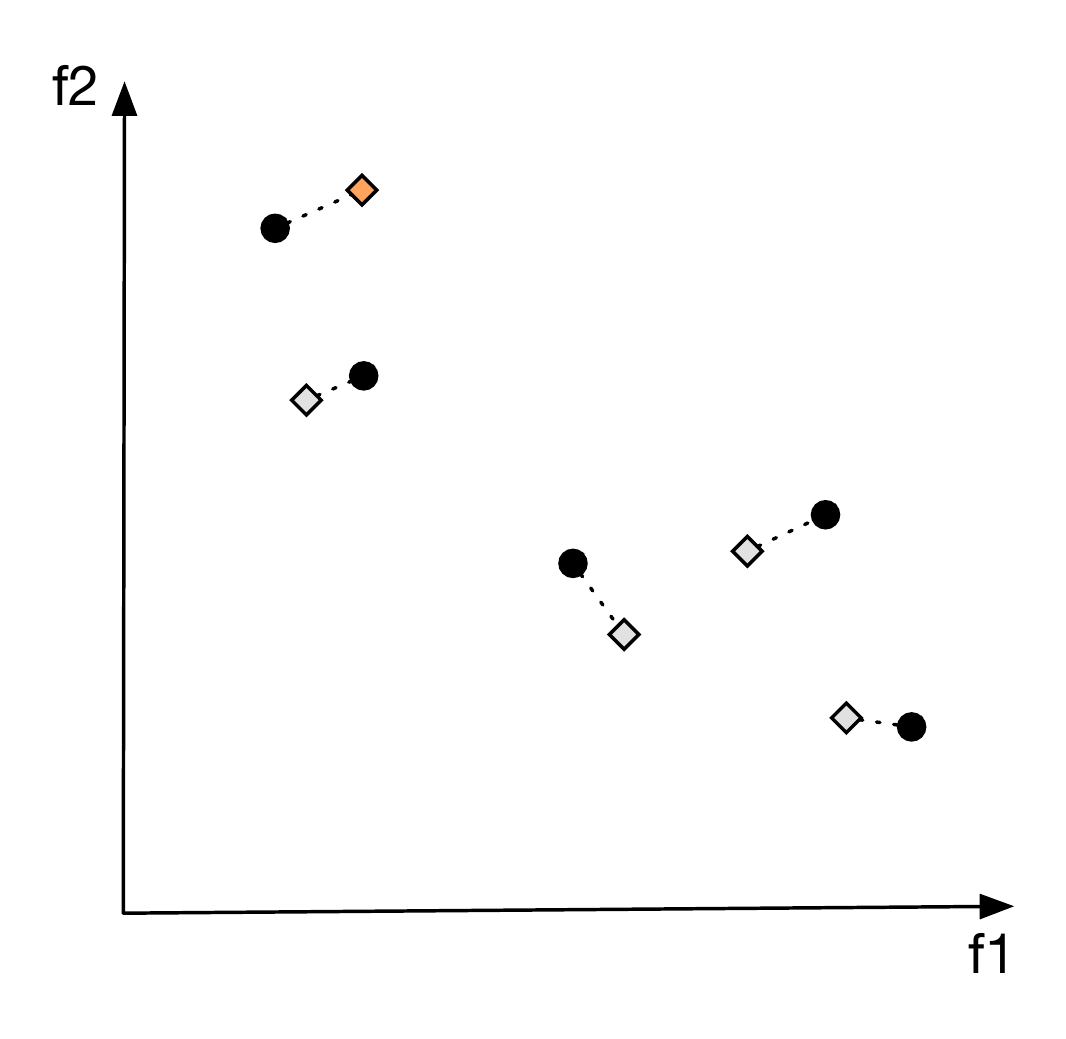}
    }\\
    \subfloat[Error by inclusion]{
    \includegraphics[width=0.3\textwidth]{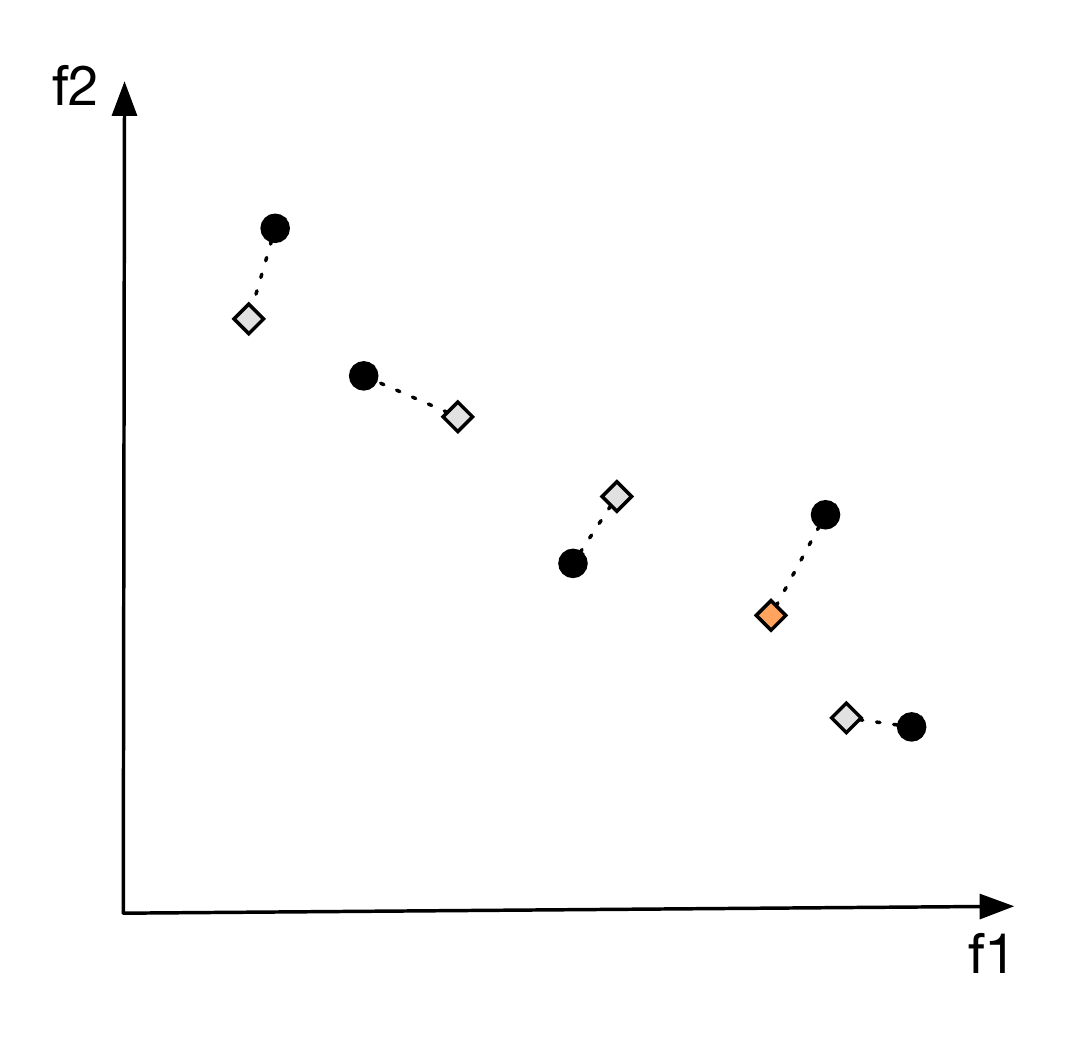}
    }\\
    \subfloat[Error by selection]{
    \includegraphics[width=0.3\textwidth]{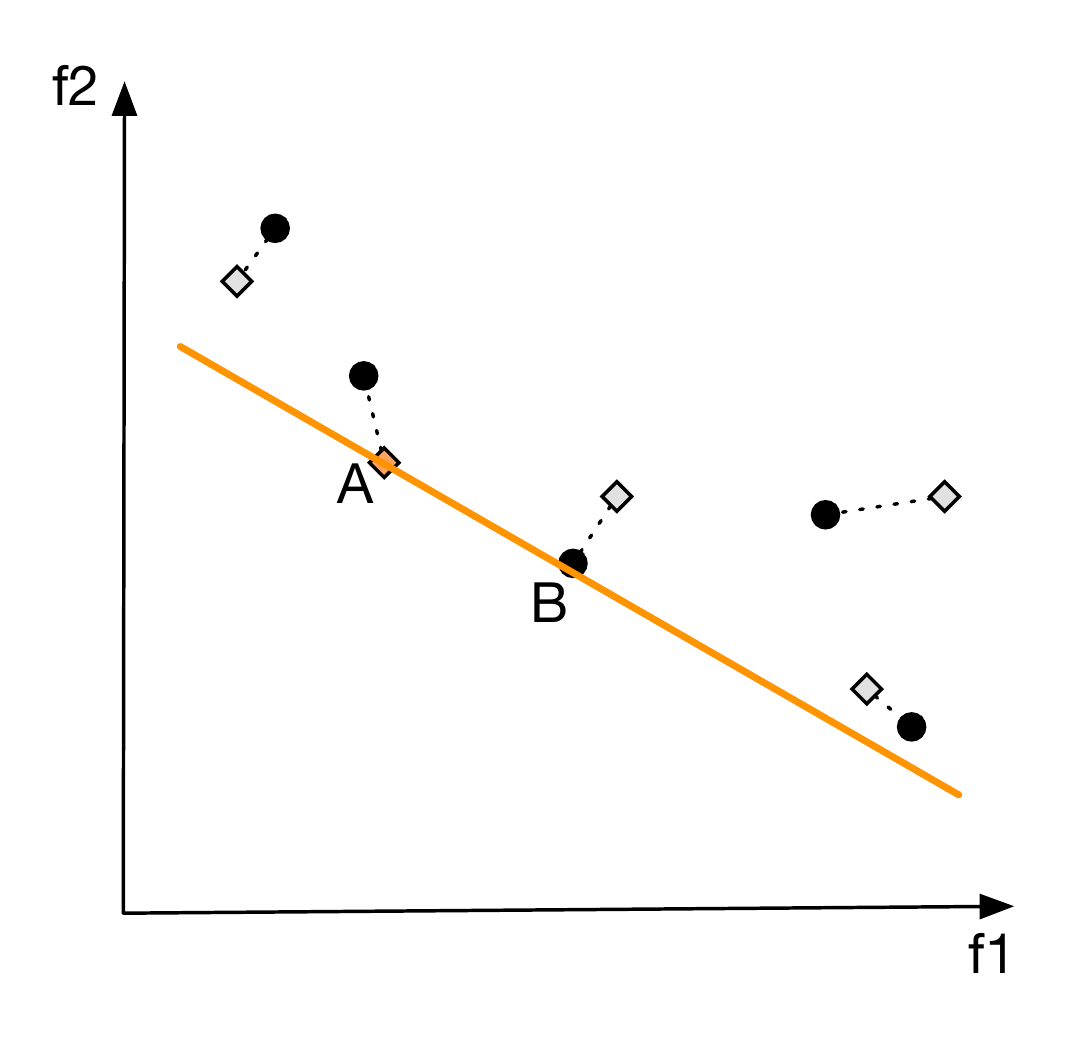}
    }
    \caption{Three possible errors due to noise. Black dots are the true objective function values of the solutions, the diamonds are the perceived objective function values. (a) The orange solution is erroneously excluded because it appears dominated while it is not.
    (b) The orange solution is erroneously included as Pareto optimal even though it is not. (c) A decision maker with a utility function represented by the organge line would erroneously select solution A even though the true utility of solution B would be higher.}
    \label{fig:errors}
\end{figure}

In multi-objective problems however, noise is more difficult to deal with.
As has already been argued by \cite{branke2007efficient}, in a setting with noise, simply using the true (noiseless) fitness values of the returned solution set would not be sufficient, because it ignores the fact that the decision maker (DM) is presented with a set of noisy and thus mis-estimated objective values for each solution, and thus there is the additional risk that the solution they pick, i.e., the solution that \emph{appears} to be the one with the highest utility, actually has a lower true utility than another solution in the set. In other words, the DM is misguided by the estimation errors and may therefore pick an inferior solution from the returned set. In fact, they might even pick a dominated solution over its dominating solution.
Thus it is clear that in a multi-objective setting under noise, the mis-estimation of the fitness values and the consequential selection errors of the DM need to be taken into account when measuring the performance of an algorithm.
This paper proposes two such metrics.

\section{Background and Literature Review}
Comparing multi-objective optimisation algorithms requires comparing sets of solutions, and various aspects are relevant, such as the closeness of the solutions to the Pareto front, the range of solutions found along each objective, and the uniformity of spread of the solutions. Various performance measures have been proposed, with the most prominent being the Hypervolume (HV) \cite{zitzler2003performance}, the Inverse Generational Distance ($IGD$) \cite{coello2004study} with its enhanced version $IGD+$
\cite{ishibuchi2015modified} and R2 \cite{hansen1994evaluating}.
A metric is called Pareto-compliant, if, whenever an approximation set $A$ dominates another approximation set $B$, the metric yields a strictly better quality value for the former than for the latter set. So far, only two performance metrics are known that are Parteto compliant, HV and $IGD^+$.

While there are many papers on multi-objective optimisation under noise (e.g., \cite{syberfeldt2010evolutionary,goh2007investigation,fieldsend2014rolling,salomon2016toolkit,teich2001pareto}), there is much less discussion on performance metrics for such problems.

Most papers in the area of noisy multi-objective optimisation just use standard performance metrics such as HV or IGD but apply them to the true fitness values of the returned set of solutions, see, e.g., \cite{chen2015performance,eskandari2009evolutionary,fieldsend2014rolling}. 

Hunter et al.~\cite{hunter2019introduction} argue that there are two errors in noisy multi-objective problems: error by inclusion (a dominated solution is falsely included in the Pareto front), error by exclusion (a non-dominated solution is falsely excluded and not returned to the DM). But even if all non-dominated solutions have been correctly identified, the DM may pick a less preferred solution because they can only make their decision based on the returned noisy objective values and not based on the true objective values. See Figure~\ref{fig:errors}

\cite{chen2015performance} propose to measure the percentage of returned solutions that are still non-dominated based on their true (undisturbed) objective values. \cite{fieldsend2014rolling} additionally propose to measure the "noise misinformation", which they define as the average distance between the returned solutions' predicted and true fitness values.
In the context of multi-objective ranking and selection, \cite{branke2019identifying} proposed to optimise a metric called hypervolume difference, motivated by the challenge for a DM to pick the correct solution in the end.


\section{Suggested Performance Metrics}
We propose two alternative performance measures, one based on IGD, the other based on the R2 metric. The basic idea is that it is necessary to simulate the DM picking a solution based on the solution set and estimated function values returned by the algorithm, then evaluating the true (noiseless) utility of this solution.

Let us assume the algorithm returns a solution set $S=\{s_1,\ldots,s_n\}$ which have true objective values $T=\{t_1,\ldots,t_n\}$ and estimated objective values $R=\{r_1,\ldots,r_n\}$, where $t_i=(f^1(s_i),\ldots f^D(s_i))$ and $r_i=(\hat f^1(s_i),\ldots \hat f^D(s_i))$.

\subsection{R2 Metric for Noisy Problems}
The R2 metric requires a parameterised utility function $U(f(x),\lambda$, and a probability distribution for the parameter $\lambda$. The simplest form is a linear combination of objectives, i.e., $U(f(x),\lambda)=\sum_{i=1}^m\lambda_i f^i(x)$. However, the concept applies equally to other utility functions such as the Tchebycheff utility of the Cobb-Douglas utility.

To calculate the R2 metric for a minimisation problem, for a set of sampled utility functions with parameters $\Lambda=\{\lambda_1,\ldots,\lambda_m\}$, the $R2$ metric
can be calculated as
\begin{eqnarray}
    R2(T,\Lambda)=\frac{1}{m}\sum_{i=1}^m \min_j\{U(t_j,\lambda_i)\},
\end{eqnarray}
i.e., for every sampled utility function, we count in the sum the utility of the best solution.

To adapt this to the noisy case, we propose the following new metric
\begin{eqnarray}
    nR2(R,\Lambda)=\frac{1}{m}\sum_{i=1}^m U(t_{j(i)},\lambda_i)
\end{eqnarray}
where $j(i)$ is the index of the solution $s_j$ with the smallest perceived utility based on the estimated fitness values, i.e.,
\begin{eqnarray}
    j(i)=\argmin_j U(r_j,\lambda_i).
\end{eqnarray}

Note that for 2 objectives and linear utility functions, the R2 metric can be calculated analytically for a continuous probability distribution for $\lambda$, rather than approximating it via Monte Carlo sampling over $\lambda$.

\subsection{Inverted Generational Distance for Noisy Problems}
The $IGD$ measure uses an approximation of the target Pareto front $A=\{a_1,\ldots,a_m\}$ of the Pareto front and then, for each solution $a\in A$ calculates the (usually Euclidean) distance to the closest point in the solution set. 
This can be put into the proposed framework where different DMs are represented the solutions in $A$. For each $a\in A$, the DM would pick the closest solution according to the estimated function values. Then, the distance of this solution's true fitness values to the desired solution $a$ is used in the calculation of $IGD$.

Let $d(a,r)$ denote the distance in objective space between two points $a$ and $r$. Then, the proposed $IGD^n$ metric of a returned solution set $S$ is defined as follows:
\begin{eqnarray}
nIGD(R,A)=\sum_{i=1}^m d(a_i,t_{j(i)})
\end{eqnarray}
where $j(i)$ is the index of the solution $s$ whose estimated function values are closest to point $a_i$, i.e., 
\begin{eqnarray}
    j(i)=\argmin_j d(a_i,r_j).
\end{eqnarray}

It is straightforward to adapt this from IGD to $IGD^+$ \cite{ishibuchi2015modified}, simply by adjusting the distance calculation and replacing  $d(a_i,t_{j(i)})$ by
\begin{eqnarray}
d_{IGD^+}(a_i, t_{j(i)})= \sqrt{\sum_{k=1}^D (\max\{t^k_{j(i)}-a_i^k,0\})^2}
\end{eqnarray}

\section{Example\label{sec:example}}
In this section, we go through a simple example, demonstrating why the standard $IGD^+$ and R2 metrics are not suitable and explaining the difference to the proposed new metrics $nIGD^+$n and $nR2$.

\begin{figure}
    \centering
    \includegraphics[width=0.4\textwidth]{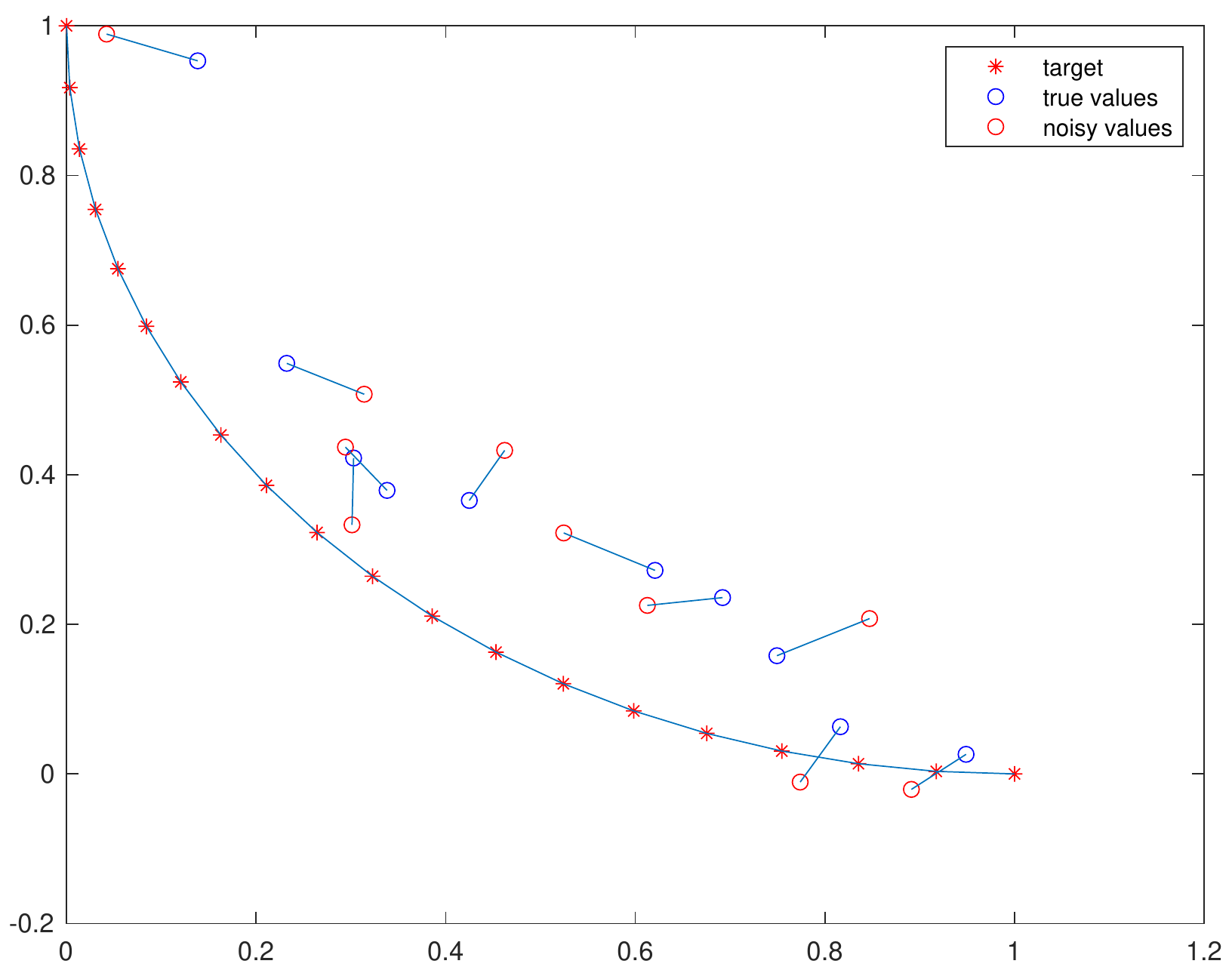}
    \caption{Observed and true quality of a set of solutions, together with true Pareto front and target values for $IGD$ calculations.}
    \label{fig:data}
\end{figure}

Figure~\ref{fig:data} shows the underlying data used: The true Pareto front (blue curve), the true quality of the solution set (blue circles), the observed quality of the solution set (red circles) and the target solutions as used by $IGD^+$ and $nIGD^+$. The true and observed qualities of each solution are connected by a line. As can be seen, the two observations at the lower end appear to be better than the true Pareto front.

\begin{figure}
    \centering
    \subfloat[R2 metric]{
    \includegraphics[width=0.4\textwidth]{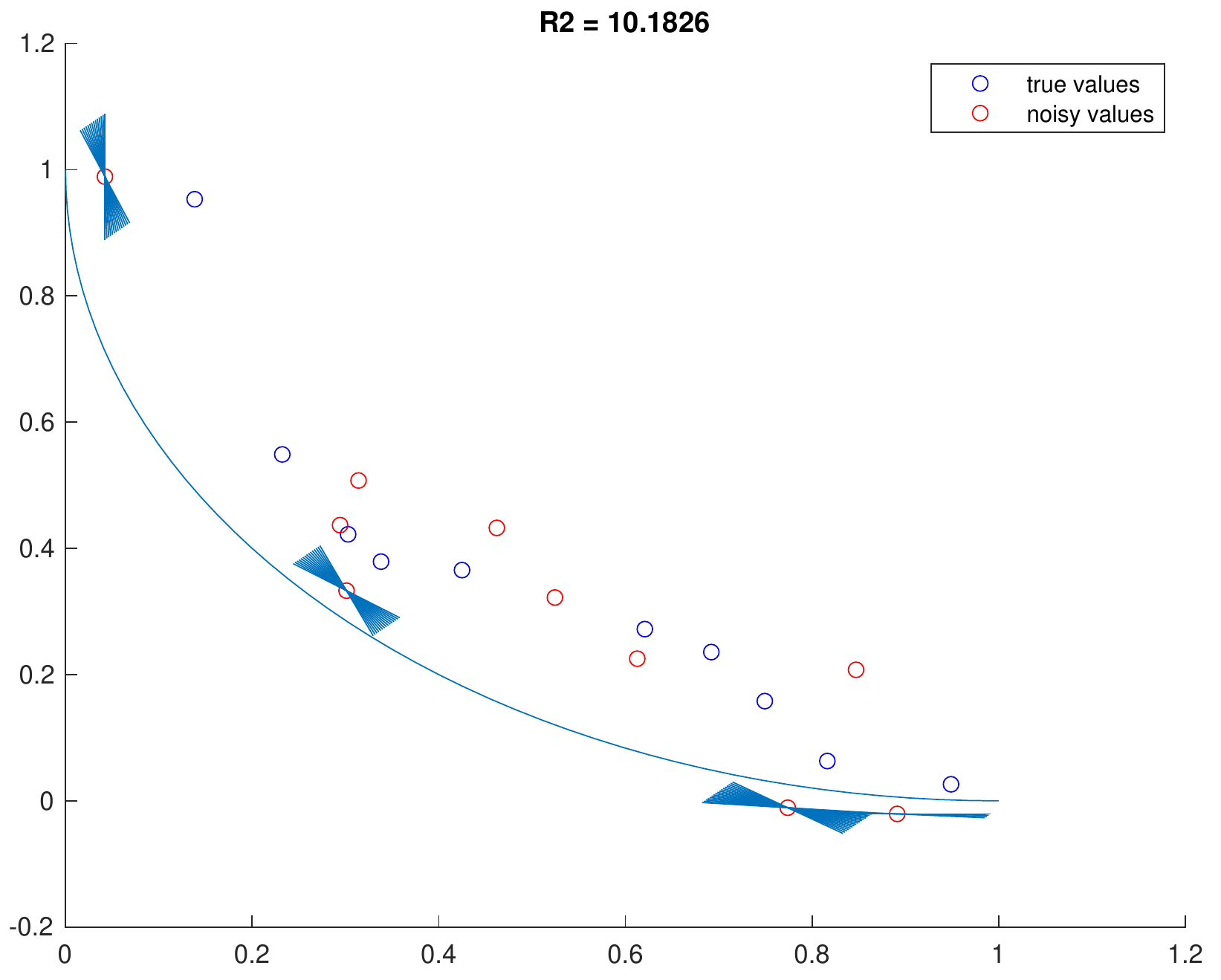}
    }\\
    \subfloat[nR2 metric]{
    \includegraphics[width=0.4\textwidth]{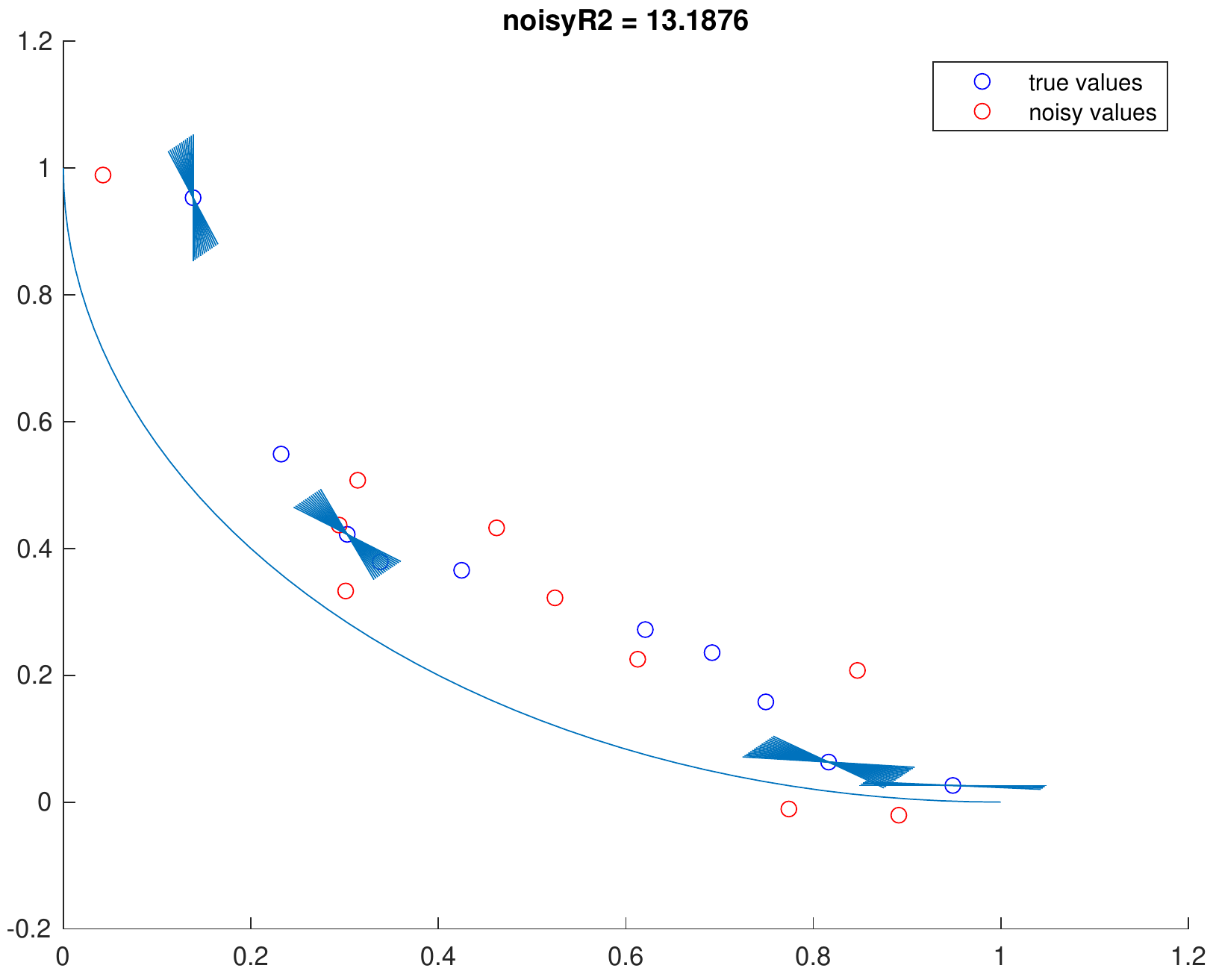}
    }
    \caption{Visualisation of (a) the R2 metric and (b) the nR2 metric.}
    \label{fig:R2}
\end{figure}

Figure~\ref{fig:R2} (a) visualises how the standard R2 metric would be calculated, using the observed solution qualities. Only four of the observed solutions are on the convex hull and will be taken into account. The short blue lines indicate the slopes of the utility functions that would favour this particular observed solution. The nR2 metric instead  uses the utilities of the true values of these very same solutions, see Figure~\ref{fig:R2} (b). Even though solution "A" is clearly non-dominated with respect to the true solution qualities, it has not been observed as such, and thus is not part of the calculation.

\begin{figure}
    \centering
    \subfloat[R2 metric Chebycheff]{
    \includegraphics[width=0.4\textwidth]{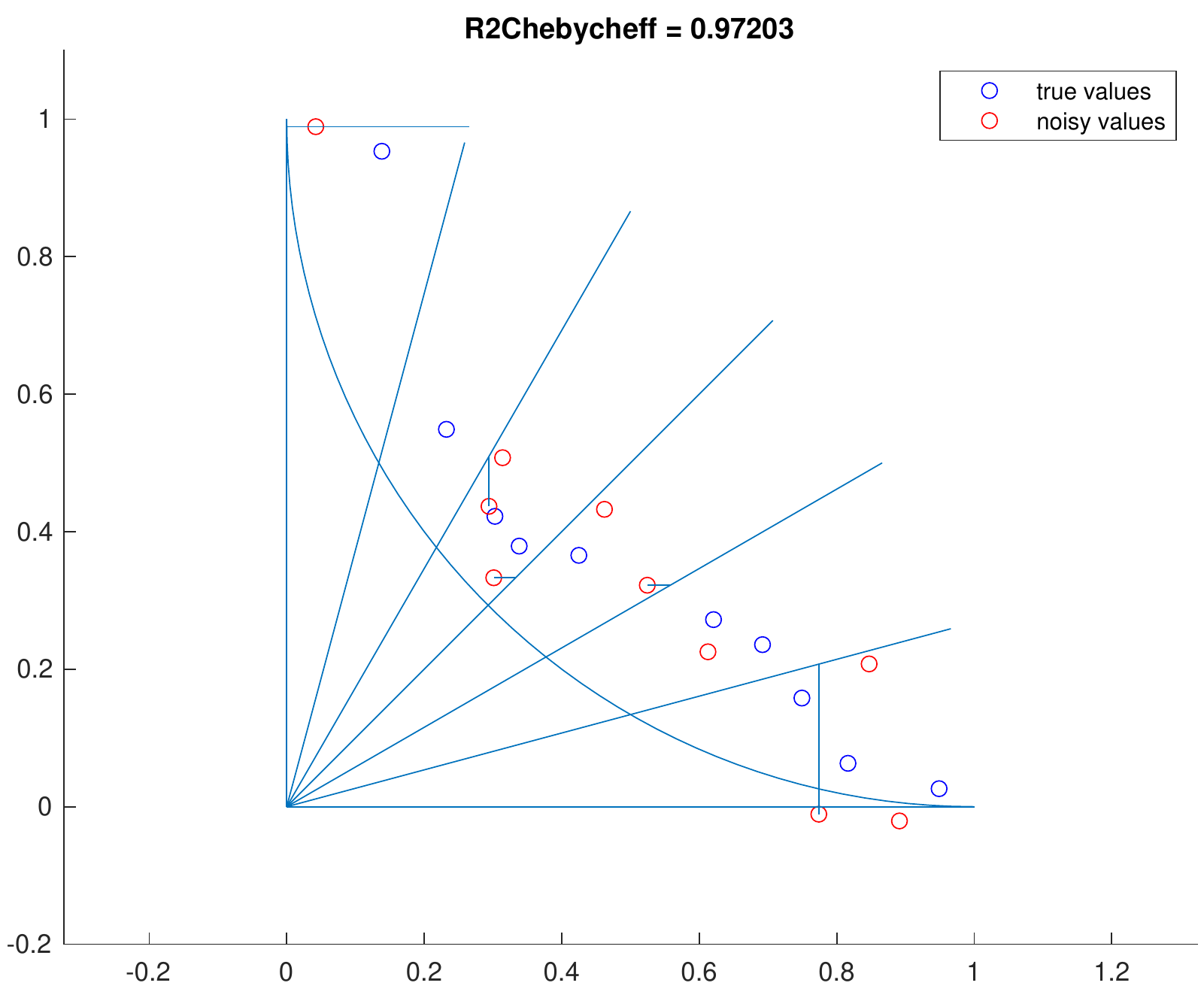}
    }\\
    \subfloat[nR2 metric Chebycheff]{
    \includegraphics[width=0.4\textwidth]{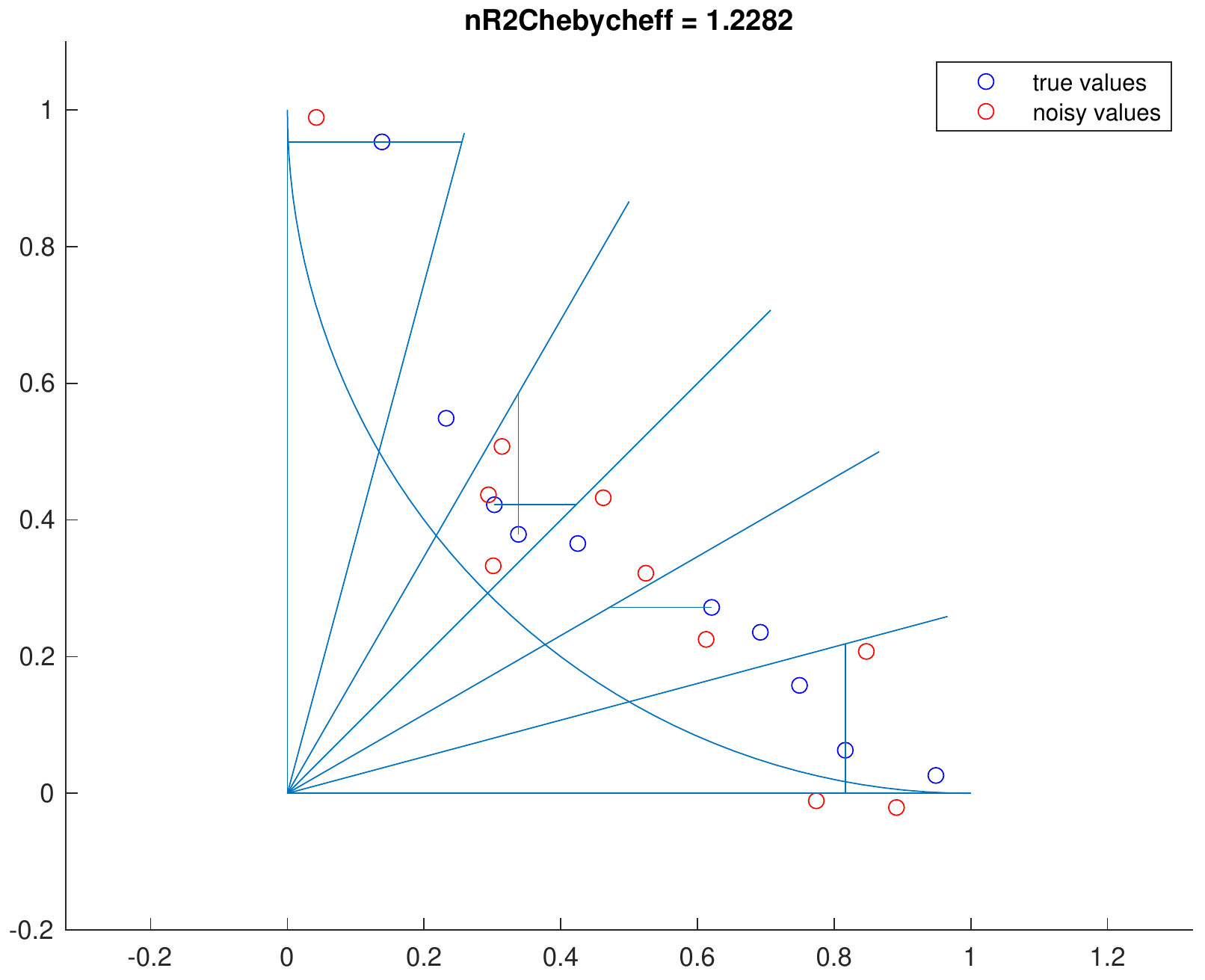}
    }
    \caption{Visualisation of (a) the R2 metric and (b) the nR2 metric when using Chebycheff scalarisation.}
    \label{fig:chebycheff}
\end{figure}

\begin{figure}
    \centering
    \subfloat[$IGD^+$ metric]{
    \includegraphics[width=0.4\textwidth]{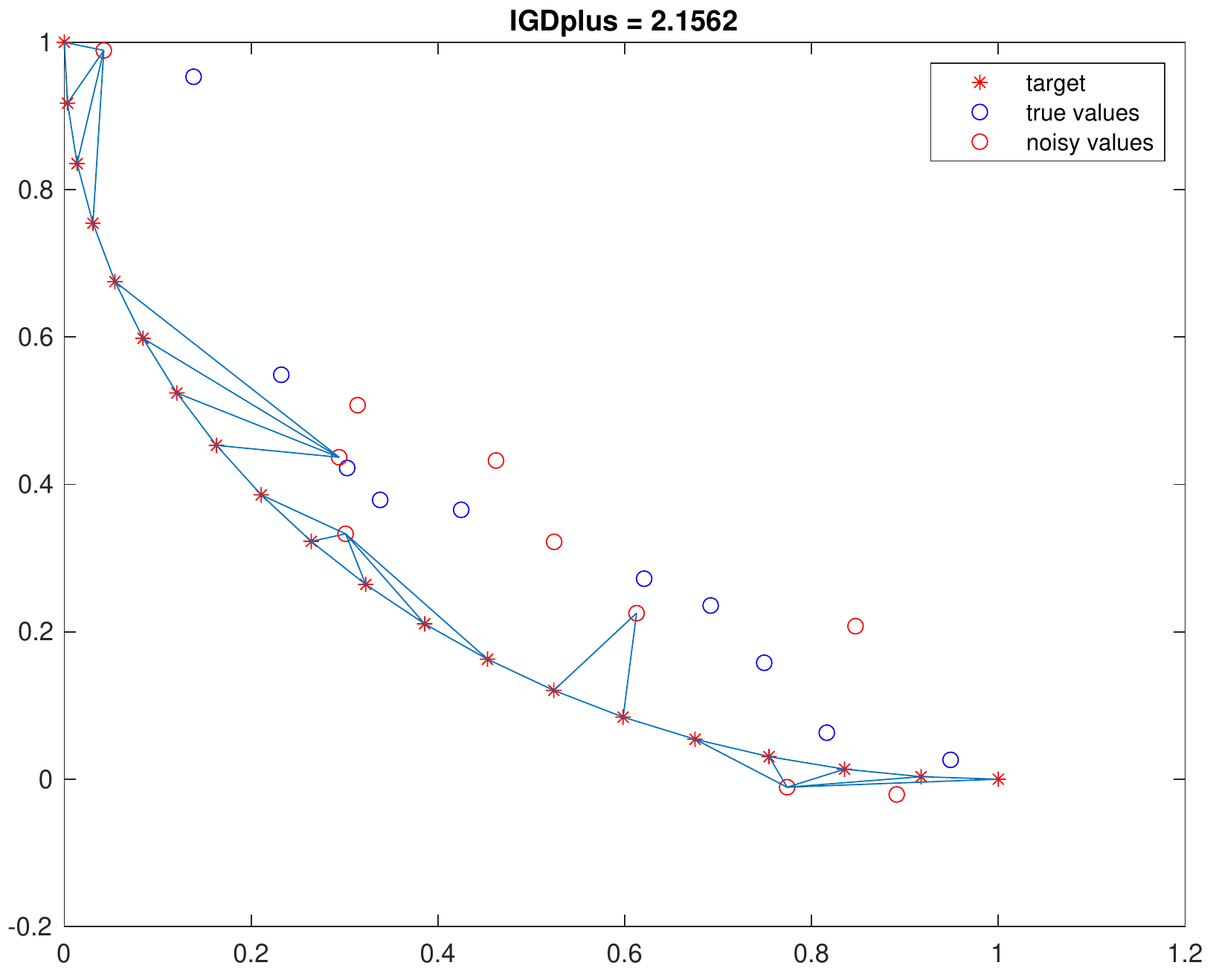}
    }\\
    \subfloat[$nIGD^+$ metric]{
    \includegraphics[width=0.4\textwidth]{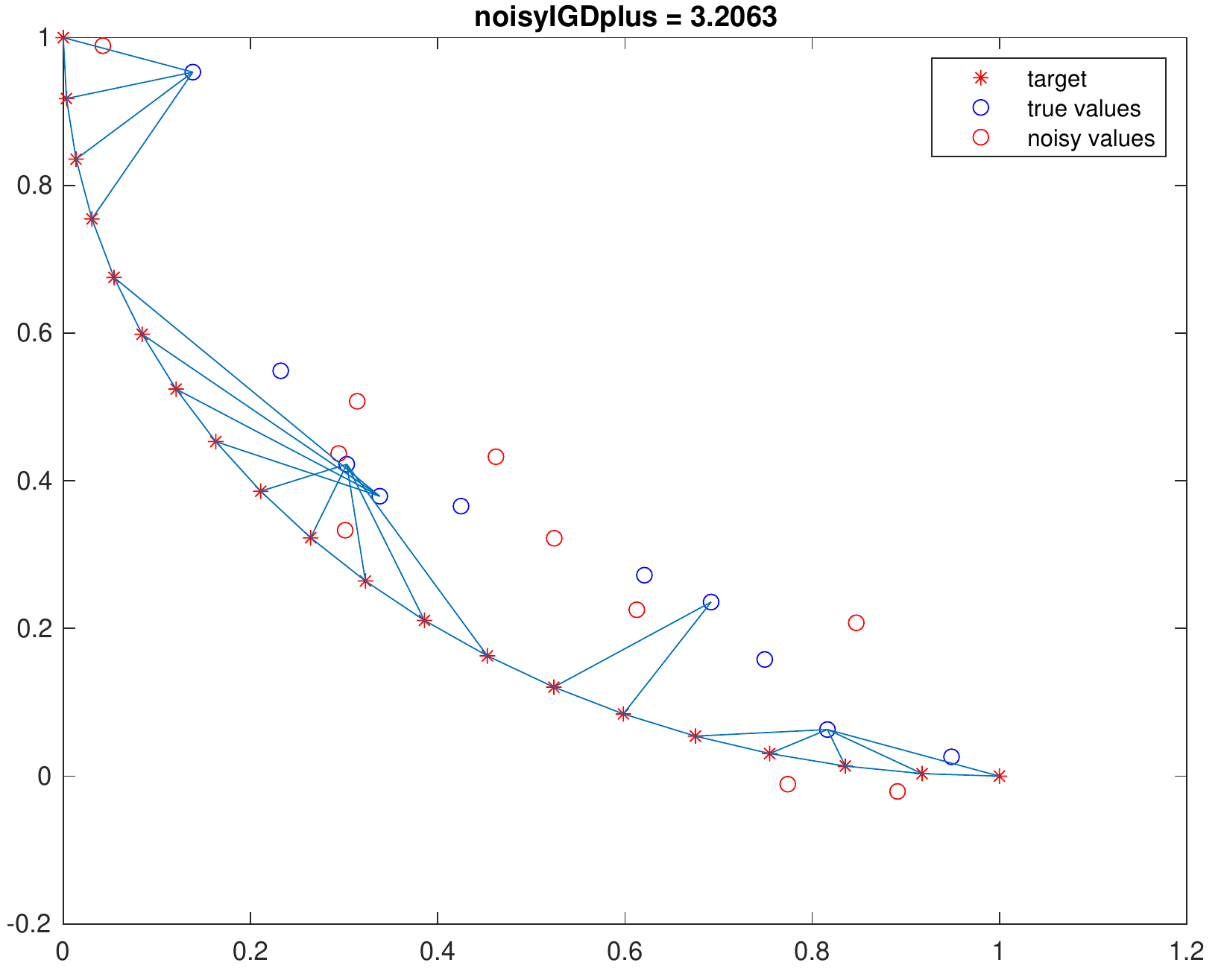}
    }
    \caption{Visualisation of (a) the $IGD^+$ metric and (b) the $nIGD^+$ metric.}
    \label{fig:igd}
\end{figure}

The $IGD^+$ metric is similarly visualised in Figure~\ref{fig:igd}. Part~(a) highlights which pairs of target and observation are taken into account for the calculation, and Part~(b) shows the same for target and true quality values as in the proposed $nIGD^+$ metric.
Where the target solution is connected to not the nearest true solution quality (which is clearly seen by lines crossing), the DM would pick the wrong solution based on the (noisy) observations.

When using Chebycheff scalarisation, the results are visualised in Figure~\ref{fig:chebycheff}. The rays represent the different weightings for the Chebycheff scalarisation, and each is connected to the solution that the DM would have picked.

Intuitively, an algorithm that is able to reduce the noise and provide more accurate predictions of a solution's true quality should be considered superior to an algorithm that reports back observed solution qualities that are far from their true qualities. 
To demonstrate that the new proposed metrics have this property, we run a test where the true solution qualities are disturbed by uniformly distributed noise in $[-\eta,\eta]^2$. The examples above show a case for $\eta=0.1$.
Table~\ref{tab:noiseImpact} reports on the computed metrics for different $\eta$ values. While the standard R2 and $IGD^+$ metric reward noise as some of the solutions appear more favourable than they are, the new nR2 and $nIGD^+$ metrics appropriately evaluate the scenarios with less noise as better.
As the noise is approaches zero, the two metrics become identical.
For the case of just using the true objective function values as is often done in the literature, the noise doesn't have any impact, so it doesn't pay off for an algorithm to improve the accuracy of the solutions' predicted objective values.

\begin{table}[h]
\caption{Different metrics tested with different levels of noise $\eta$. Lower values are better in all cases. R2 is the R2 metric with linear utility functions, R2c is the R2 metric with Chebycheff utility function.\label{tab:noiseImpact}}
    \centering
    \begin{tabular}{c|cccccc}
         $\eta$&R2&nR2&R2c&nR2c&$IGD^+$&$nIGD^+$  \\\hline
         0.01&12.7475 & 12.9576 & 1.1083 &1.1148& 6.4231 & 6.4411\\
         0.05&11.6563 & 13.0960 & 1.0562 & 1.1689 & 6.1052 & 6.8660\\
         0.1&10.1826 &13.1876 & 0.9720 & 1.2282 & 5.5193 & 8.1651\\
         0.2& 7.0602 & 13.4119 & 0.8372 & 1.3339 & 4.2270 & 8.3205\\\hline
    \end{tabular}
    
\end{table}

\section{Discussion}
As explained above, in noisy multi-objective optiimisation it is important to take into account the selection error of the DM due to wrongly predicted objective function values. This is straightforward for the R2 metric, because this metric samples from DM utility functions and thus the selection process of the DM can be directly modelled. A pre-requisite for this metric is however that we can identify a probability distribution over the parameters of a representative utility function model. Using a linear utility function also only rewards solutions that are predicted to be on the convex hull.

It is also possible based on the $IGD$ metric, if one assumes that the target approximation of the Pareto optimal front represents different possible DM preferences, and a DM with a particular target solution would select the closest predicted solution. While the DM model (picking the closest predicted solution) seems straightforward and natural, the $IGD$ metric requires an approximation of the true Pareto frontier, which obviously is not available for real-world problems.

It is less obvious but interesting for future work to also find an adaptation of the popular HV metric.

Note that the metrics proposed can compare solution sets of different sizes, so an interesting question for an algorithm designer is how many solutions perceived as Pareto-optimal the algorithm should return. It doesn't make sense to return a solution that is dominated by another solution in the returned set, as this dominated solution would never enter the calculation. However, it can be beneficial to for example remove even predicted non-dominated solutions from the set if their predicted objective function values are very uncertain, as their true values may turn out much worse than predicted, deteriorating the performance metric.

\section{Conclusion}
We have pointed out that when evaluating an algorithm for noisy multi-objective optimisation, one should take into account the utility loss for a DM who may pick an inferior solution because the algorithm didn't accurately predict the quality of the identified solutions. Based on this observation, we proposed modifications of the well-know R2 and $IGD/IGD^+$ metrics and demonstrated with an example that they have the desired property.
As future work, it may be worthwhile to also develop an adaptation of the popular hypervolume metric to the noisy case. Also, it would be helpful to have a performance metric that doesn't rely on knowing the solutions' true objective function values.


\bibliographystyle{IEEEtran}
\bibliography{literature}

\end{document}